\documentclass[12pt]{amsart}
\usepackage{amsfonts,latexsym} 
\usepackage{url}
\usepackage{graphicx,graphics,hieroglf,skak}
\usepackage{amsthm}
\usepackage{amssymb, amsmath}
\usepackage{url,enumitem}
\usepackage[cm]{fullpage}

\theoremstyle{definition}

\theoremstyle{remark}

\begin{document}

	\begin{center} 

	{\bf RAMANUJAN, ROBIN, HIGHLY COMPOSITE NUMBERS,\\AND THE RIEMANN HYPOTHESIS}\\
\vspace{5 mm}

Jean-Louis Nicolas\\
Universit\'e de Lyon; CNRS; Universit\'e Lyon 1;\\
Institut Camille Jordan, Math\'ematiques,\\
21 Avenue Claude Bernard, F-69622 Villeurbanne cedex, France\\
E-mail: nicolas@math.univ-lyon1.fr\\
\vspace{5 mm}

Jonathan Sondow\\
209 West 97th Street \#6F,
New York, NY 10025, USA\\
E-mail: jsondow@alumni.princeton.edu\\

		\end{center}

\vspace{5 mm}



\noindent 
{\Small {\bf Abstract.} We provide an historical account of equivalent conditions for the Riemann Hypothesis arising from the work of Ramanujan and,
later, Guy Robin on generalized highly composite numbers. The first part of the paper is on the mathematical background of our subject. The
second part is on its history, which includes several surprises.}\\
\bigskip

\begin{center}
{\bf 1. MATHEMATICAL BACKGROUND}\\
\end{center}
\bigskip

\bigskip

\noindent{\bf Definition.} The \emph{sum-of-divisors function} $\sigma$ is defined by
\begin{center}
$\displaystyle\sigma(n):=\sum_{d \mid n}d = n\sum_{d \mid n}\frac1d.$
\end{center}

In 1913, Gr\"{o}nwall found the maximal order of~$\sigma$.\\

\noindent{\bf Gr\"{o}nwall's Theorem \cite{GRO}.} {\it The function}
\begin{center}
$G(n):=\dfrac{\sigma(n)}{n \log\log n} \qquad(n>1)$
\end{center}
{\it satisfies}
\begin{center}
$\displaystyle\limsup_{n\to\infty}\, G(n) =  e^\gamma = 1.78107\dotso,$
\end{center}
{\it where$$\gamma:=\lim_{n\to\infty} \left(1+\frac12+\dotsb+\frac1n-\log n\right)=0.57721\dotso$$is the Euler-Mascheroni constant.}\\

\centerline{\includegraphics[width=1.8in]{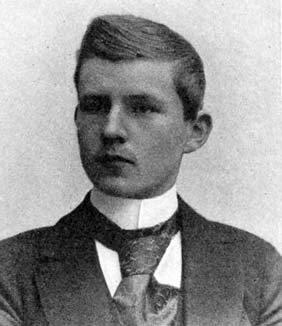}}\begin{center}{ Thomas Hakon GR\"{O}NWALL (1877--1932)}\end{center}

\newpage

Gr\"{o}nwall's proof uses:\\

\noindent{\bf Mertens's Theorem \cite{MER}.} 
{\it If $p$~denotes a prime, then}
$$\lim_{x\to\infty} \frac{1}{{\log x}} \prod_{p\le x}\left(1-\frac1p\right)^{-1}= e^{\gamma}.$$

\centerline{\includegraphics[width=2in]{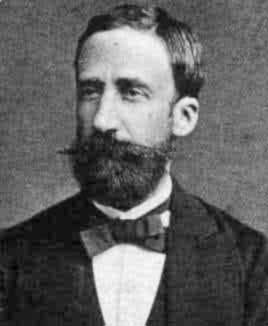}}\begin{center}{ Franz MERTENS (1840--1927)}\end{center}

\bigskip

\noindent{\bf Ramanujan's Theorem \cite{RAMNarosa, RAM97, BerAnd}.}
{\it If RH is true, then for $n_0$ large enough,}
\begin{equation}\label{GRAM}
n>n_0 \implies G(n) < e^\gamma.\\
\end{equation}

\centerline{\includegraphics[width=3in]{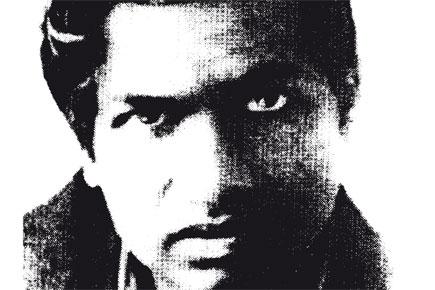}}\begin{center}{Srinivasa RAMANUJAN (1887--1920)}\end{center}

\bigskip

To prove that RH implies \eqref{GRAM}, Ramanujan introduces a real non-negative
parameter $s$, considers the multiplicative function $n
\mapsto\sigma_{-s}(n)=\sum_{d\mid n} d^{-s}$ (for example, $\sigma_{-1}(n)=\sigma(n)/n$),  
and calls an integer $N$ a {\it generalized highly composite number} if
$$N' < N \implies \sigma_{-s} (N') < \sigma_{-s} (N).$$
When $s=1$ these numbers have been called {\it superabundant} by Erd\H
{o}s and Alaoglu,  while for  $s\neq 1$ they have only
been studied by Ramanujan.
Further, Ramanujan calls an integer $N$ a {\it generalized superior
  highly composite number} of parameter $\varepsilon >0$ if
$$N' < N \implies\frac{\sigma_{-s}(N)}{N^\varepsilon} \geq
\frac{\sigma_{-s}(N')}{(N')^\varepsilon} $$
and
$$N' > N \implies\frac{\sigma_{-s}(N)}{N^\varepsilon} >
\frac{\sigma_{-s}(N')}{(N')^\varepsilon} .$$
When $s=1$
these numbers have been called {\it colossally abundant} by Erd\H
{o}s and Alaoglu. It is easily seen that all generalized superior
highly composite numbers are generalized highly composite.

The prime factorization of a generalized superior
highly composite number $N$ can be obtained from the value of the
parameter $\varepsilon$. For $r=1,2,3\ldots$, Ramanujan defines $x_r$
by
$$x_r^\varepsilon=\frac{1-x_r^{-s(r+1)}}{1-x_r^{-sr}}$$
and then 
$$N=\prod_{r=1}^R e^{\vartheta(x_r)}$$
where $\vartheta(x)=\sum_{p\leq x} \log p$ denotes Chebyshev's function
and $R$ is the largest integer such that $x_R\geq 2$. One has
$$\sigma_{-s}(N)=\prod_{r=1}^R \prod_{p\leq x_r}
\frac{1-p^{-s(r+1)}}{1-p^{-sr}},$$
and, to estimate $\sigma_{-s}(N)$, one needs an estimate of
$$\sum_{p\leq x} \log\left( 1-\frac{1}{p^s}\right)=-\sum_{p\leq x}
\int_s^\infty \frac{\log p}{p^t-1} {\rm d}t$$
whence the idea of considering the sum $\sum_{p\leq x} \frac{\log
  p}{p^s-1}\cdot$

\bigskip

Here is an excerpt from Ramanujan's proof.\\

\noindent Ramanujan \cite[p.~133]{RAM97}: {\it \dotso assume that $\dotso s>0\dotso$ if $p$ is the largest prime not greater than~$x$, then}

\begin{align*}
\frac{\log2}{2^s-1}+\frac{\log3}{3^s-1}+\frac{\log5}{5^s-1}+\dotsb+\frac{\log p}{p^s-1}
=C+\int^{\theta(x)} \negthickspace\frac{dx}{x^s-1} - s\negthickspace\int \frac{x-\vartheta(x)}{x^{1-s}(x^s-1)^2}\,dx + O\{x^{-s}(\log x)^4\}.
\end{align*}

\bigskip

\noindent {\it But it is known that}
\begin{equation*}
x-\theta(x) = \sqrt x + x^{\frac13}+ \sum \frac{x^\rho}{\rho} - \sum \frac{ x^{\frac12\rho}}{\rho} + O(x^{\frac15})
\end{equation*}

\noindent {\it where $\rho$ is a complex root of} $\zeta(s)$.\dotso\\

The last equation is a variant of the classical explicit formula in prime number theory. This shows ``explicitly'' how Ramanujan used RH in his proof.\\


From the estimate for $\sigma_{-s}(N)$, Ramanujan deduces that, for
$1/2< s < 1$ and all sufficiently large integers $n$, the upper bound 
\begin{eqnarray*}
\sigma_{-s}(n)  &\leq & |\zeta(s)| \exp\left\{{\rm Li}((\log
n)^{1-s})-\frac{2s(2^{\frac{1}{2s}}-1)}{2s-1} \frac{(\log n)^{\frac12-s}}{\log
\log n}\right\}\\
& & -\frac{s}{\log \log n}\sum_\rho \frac{(\log
  n)^{\rho-s}}{\rho(\rho-s)}+O\left\{\frac{(\log n)^{\frac{1}{2}-s}}{(\log \log n)^2}\right\}
\end{eqnarray*}
holds. Finally, making $s$ tend to $1$, he gets
$$\limsup_{n\to\infty} \, (\sigma_{-1}(n)-e^\gamma \log \log n)\sqrt{\log n}\leq
-e^\gamma(2\sqrt 2-4-\gamma+\log4\pi)=-1.393\ldots <0.$$
Since $G(n)=\sigma_{-1}(n)/\log\log n$, this proves \eqref{GRAM}.


\bigskip

\noindent{\bf Robin's Theorem \cite{RobDPP,RobJMPA}.}
{\it RH is true if and only if}
\begin{equation*}\label{GROB}
n > 5040\ (= 7!)\implies G(n) < e^\gamma.
\end{equation*}

\bigskip

\centerline{\includegraphics[width=2in]{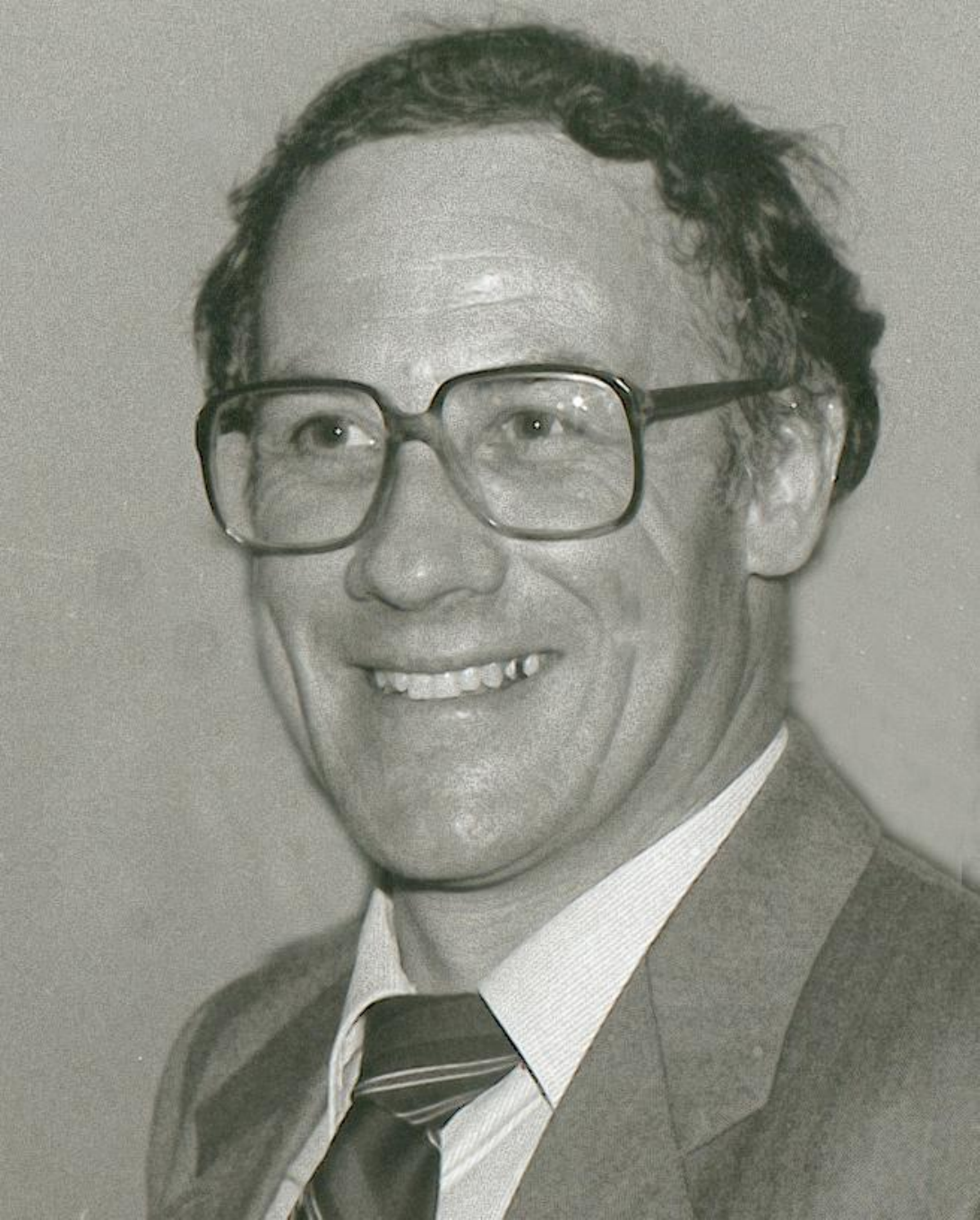}}\begin{center}{ Guy ROBIN}\end{center}

\bigskip

\centerline{\includegraphics[width=4in]{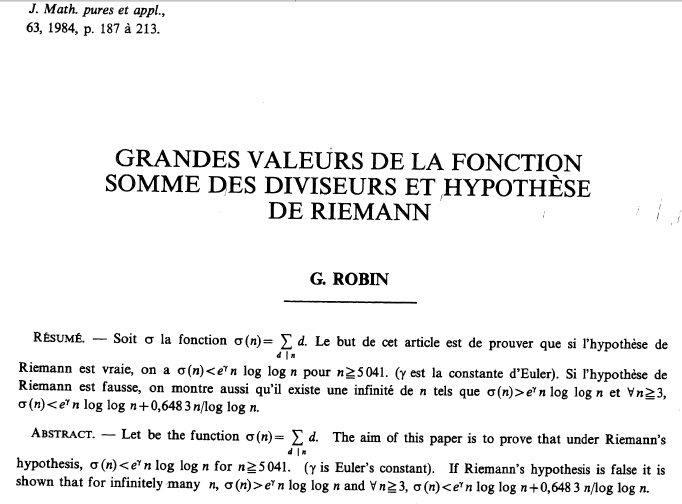}}
\begin{center}{ Robin's paper on $\sigma$ and RH, Journal de Math\'{e}matiques Pures et Appliqu\'{e}es, 1984}\end{center}
\vspace{.3in}


To prove his theorem, Robin uses generalized superior highly
composite numbers only with $s=1$, i.e., colossally abundant
numbers (CA for short). First, he shows that if $N' < N''$ are two
consecutive CA numbers, then
$$N' < n < N'' \implies G(n) \leq \max(G(N'),G(N'')).$$
Second, by numerical computation, he checks that $G(N) <
e^\gamma$ for all integers $N$ with
$5041 < N < 55440$, as well as for all CA numbers $N\geq
55440$ whose largest prime factor $P^+(N)$ is $< 20000$.

Further, if a CA number $N$ satisfies $P^+(N) > 20000$, then getting an
upper bound for $\sigma_{-1}(N)$ requires a precise estimate of the
Mertens product
\begin{equation}
  \label{MERT}
\prod_{p\leq x} \left(1-\frac1p\right)^{-1}.  
\end{equation}

\newpage

The sum-of-divisors function $\sigma$ and Euler's totient function~$\phi$, defined as
\begin{equation*}
\displaystyle\phi(n):=\sum_{\stackrel{1\le k\le n}{(k,n)=1}}1= n\prod_{p \mid n}\left(1-\frac1p\right),
\end{equation*}
are related by the inequalities
\begin{equation*}
\frac{6}{\pi^2} < \frac{\sigma(n)}{n}\cdot\frac{\phi(n)}{n} < 1,
\end{equation*}
which hold for all $n>1$.
Mertens's Theorem implies that the minimal order of $\phi$ is given by
$$\displaystyle\limsup_{n\to\infty}\, \frac{n/ \log\log n}{\phi(n)} =  e^{\gamma}.$$

To estimate the product \eqref{MERT},
Robin used ideas from a result on the $\phi$~function proved by his thesis advisor Nicolas in 1983.\\

\noindent{\bf Nicolas's Theorem \cite{NicDPP,NicJNT}.}
{\it RH is true if and only if}
$$prime\ p>2\implies \frac{p\#/\log\log p\#}{\phi(p\#)} > e^{\gamma}, $$
{\it where $p\# := 2\cdot3\cdot5\cdot7\cdot11\dotsb p$ denotes a primorial.}

\bigskip

\centerline{\includegraphics[width=2.5in]{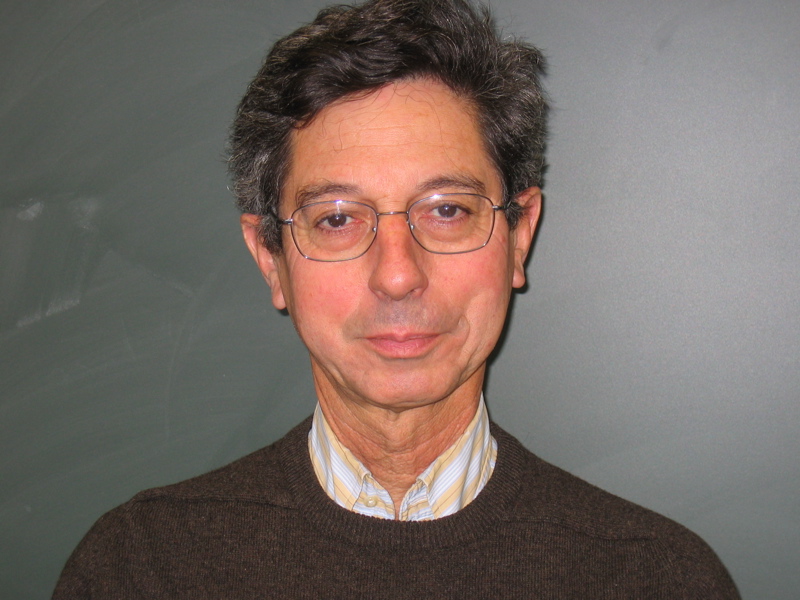}}\begin{center}{ Jean-Louis NICOLAS}\end{center}

\bigskip

Nicolas in turn used Landau's Oscillation Theorem \cite{Lan}, which Landau had applied in 1905 to prove a form of Chebyshev's bias $\pi(x;4,3)>\pi(x;4,1)$.\\

\centerline{\includegraphics[width=2in]{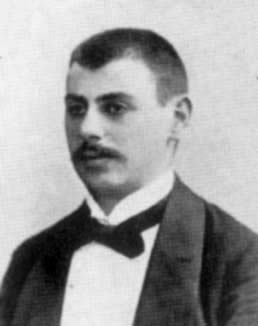}}\begin{center}{ Edmund LANDAU (1877--1938)}\end{center}

\bigskip



\noindent{\bf Caveney, Nicolas, and Sondow's Theorems \cite{CNS,CNS2}.}\\
{\it Define an integer $N>1$ to be a~\emph{GA1 number} if $N$~is composite and $G(N) \ge G(N/p)$ for all prime factors~$p$. Call an integer $N$ a \emph{GA2 number} if $G(N) \ge G(aN)$ for all multiples $aN$. Then:}\\
1. {\it RH is true if and only if $4$~is the only number that is both GA1 and GA2.}\\
\noindent2. {\it A GA2 number $N>5040$ exists if and only if RH is false, in which case $N$ is even and} $>10^{\,8576}$.


\bigskip

\centerline{\includegraphics[width=2in]{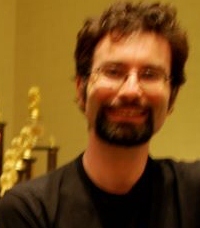}}
\begin{center}{ Geoffrey CAVENEY}\end{center}

\vspace{.5in}

\begin{center}
{\bf 2. HISTORY}
\end{center}

\bigskip
Our story begins in 1915, when Ramanujan published the first part of his dissertation ``Highly Composite Numbers'' (HCN for short).\\

\centerline{\includegraphics[width=4in]{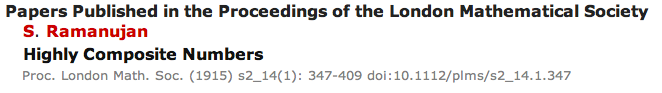}}\begin{center}{Ramanujan's HCN Part 1, Proceedings of the London Mathematical Society, 1915}\end{center}
\bigskip

\noindent Ramanujan (in \cite{RAM15}):
{\it I define a highly composite number as a number whose number of divisors exceeds that of all its predecessors.}\\

His thesis was written at Trinity College, University of Cambridge, where his advisors were Hardy and Littlewood.


\centerline{\includegraphics[width=2in]{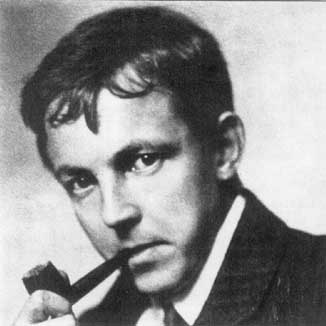}}\begin{center}{ Godfrey Harold HARDY (1877--1947)}\end{center}\

\centerline{\includegraphics[width=1.6in]{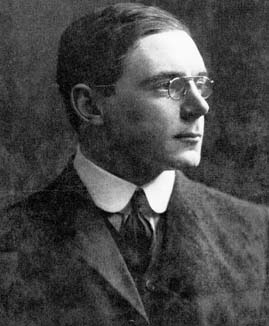}}\begin{center}{ John Edensor LITTLEWOOD (1885--1977)}\end{center}\


\bigskip

\centerline{\includegraphics[width=3.3in]{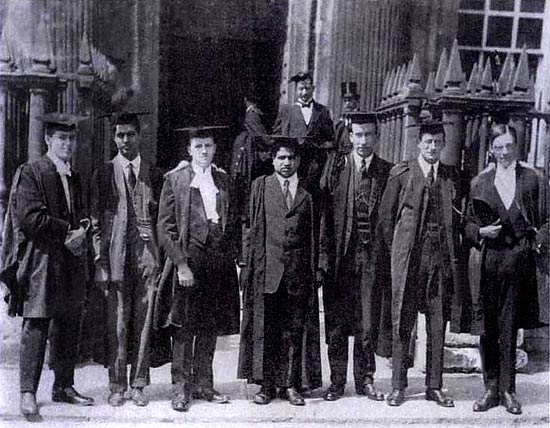}}\begin{center}{ Ramanujan (center) at his degree ceremony, 1916}\end{center}


\bigskip

In 1944, Erd\H{o}s published a paper ``On highly composite and similar numbers'' with Alaoglu \cite{AE}.\\


\noindent Erd\H{o}s (in ``Ramanujan and I'' \cite{ERD}): {\it Ramanujan wrote a long paper \cite{RAM15} on this subject. Hardy rather liked this paper but perhaps not unjustly called it nice but in the backwaters of mathematics. $\dotso$ Ramanujan had a very long manuscript on highly composite numbers but some of it was not published due to a paper shortage during the First World War.  }

\bigskip

\centerline{\includegraphics[width=1.5in]{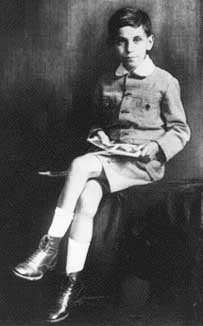}}\begin{center}{ Paul ERD\H{O}S (1913--1996)}\end{center}

\bigskip




\noindent Dyson (email to Sondow, 2012): {\it Hardy told me, ``Even Ramanujan could not make highly composite numbers interesting.'' He said it to discourage me from working on H. C. numbers myself. I~think he was right.}\\

\centerline{\includegraphics[width=1.5in]{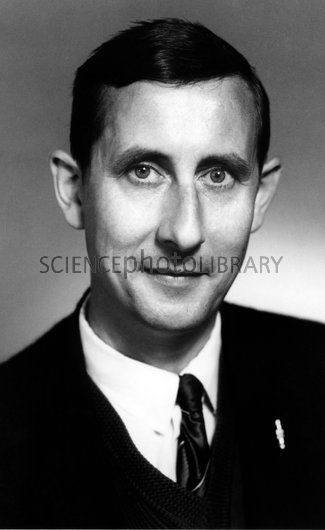}}
\vspace{-5mm}
\begin{center}{ Freeman DYSON}\end{center}




In 1982 Rankin published a paper on ``Ramanujan's manuscripts and notebooks.'' He quoted Hardy's mention of ``the suppressed part of HCN'' in a 1930 letter to Watson.\\

\noindent Rankin (in \cite{RAN}): {\it The most substantial manuscript consists of approximately 30 pages of HCN carrying on from where the published paper stops.}\\

\centerline{\includegraphics[width=2in]{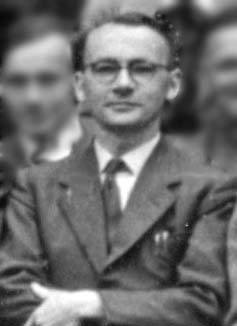}}
\vspace{-5mm}
\begin{center}{Robert Alexander RANKIN (1915--2001)}\end{center}\






By a curious coincidence, 1981--1982 is also the year of S\'{e}minaire Delange-Pisot-Poitou's exposition \cite{RobDPP} of Robin's Theorem, in which he improved on Ramanujan's Theorem 
without ever having heard of it!\\

\noindent Berndt (email to Sondow, 2012): {\it It is doubtful that Rankin took notice of Robin's paper. I definitely did not.}\\

Thus, Rankin and Berndt on the English-American side, and Nicolas and Robin on the French side, were not communicating.\\

\centerline{\includegraphics[width=2.3in]{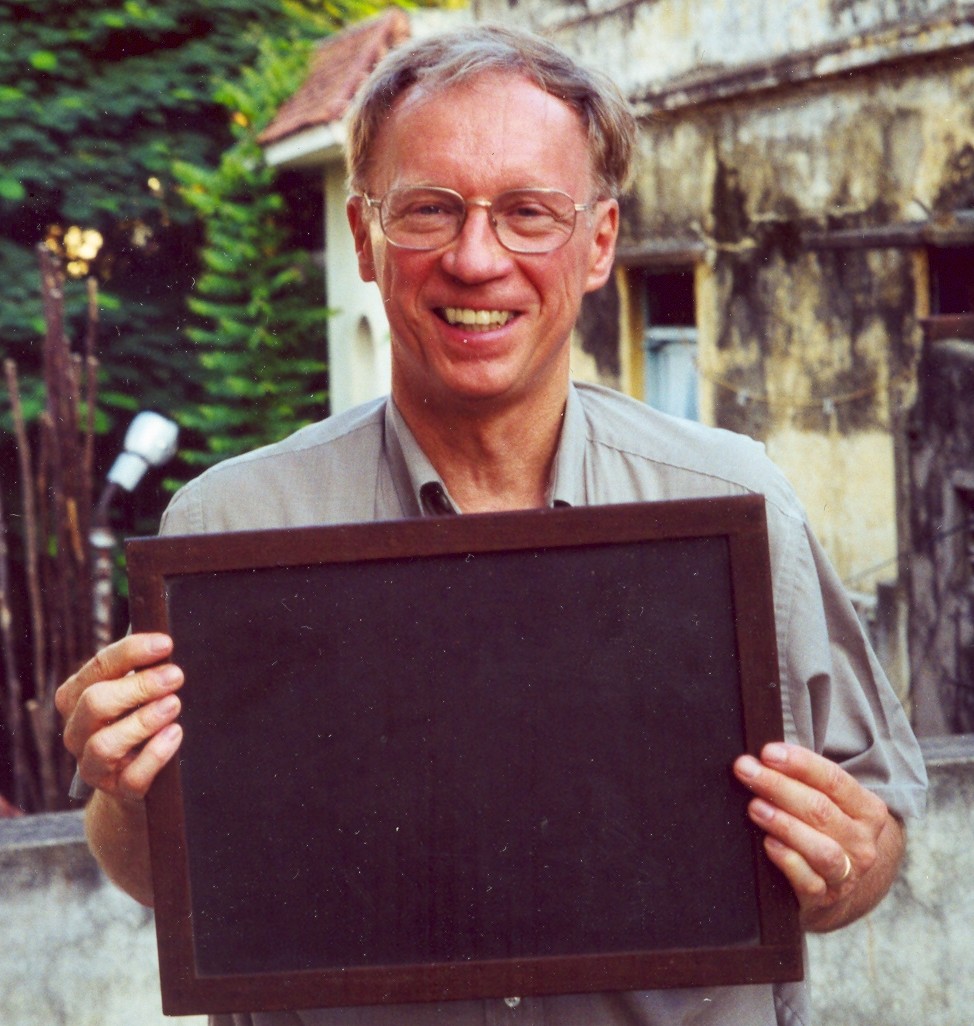}}\begin{center}{ Bruce Carl BERNDT holding Ramanujan's slate}\end{center}\

\noindent Berndt (email to Sondow, 2012):
{\it After I began to edit Ramanujan's notebooks, I wrote Trinity College  in
 1978 for a copy of the notes that Watson and Wilson made in their  efforts
 to edit the notebooks.  I also decided to write for copies of all the
 Ramanujan material that was in the Trinity College Library. Included in
 their shipment to me was the completion of Ramanujan's paper on highly
 composite numbers.  I put all of this on display during the Ramanujan
 centenary meeting at Illinois in June, 1987.}




\bigskip

\noindent Nicolas (email to Sondow, 2012):
{\it I keep a very strong souvenir of the conference organised in
 Urbana-Champaign in 1987 for the one hundred anniversary of
 Ramanujan. It is there that I discovered the hidden part  of ``Highly
 Composite Numbers''} [first published in \cite{RAMNarosa}, later in \cite{RAM97}, and again in \cite{BerAnd}]{\it.

 What I have not written is that there was an error of calculus in
 Ramanujan's manuscript which prevented him from seeing Robin's
 Theorem. Soon after discovering the hidden part, I read it and saw
 the difference between Ramanujan's result and Robin's one. Of course,
 I would have bet that the error was in Robin's paper, but after
 recalculating it several times and asking Robin to check, it turned
 out that there was an error of sign in what Ramanujan had written.}

\bigskip
 
 Thus it happened that Robin avoided the fate of many mathematicians, who have found that (Berndt \cite{BCCL}, \cite{Ber} quoting Gosper): {\it Ramanujan reaches his hand from his grave to snatch your theorems from you.}

\bigskip
 
Ramanujan's Theorem was not explicitly stated by him in HCN Parts~1 or~2. Nicolas and Robin formulated it for him in Note~71 of their annotated and corrected version of HCN Part 2.\\

\noindent Nicolas and Robin (in \cite{RAM97}): {\it It follows from (382)} [(the corrected version of Ramanujan's  formula)] {\it that under the Riemann hypothesis, and for $n_0$ large enough,}
$$n>n_0 \implies \sigma(n)/n < e^\gamma\log\log n.$$

\noindent{\it It has been shown in \cite{RobJMPA} that the above relation with $n_0=5040$ is equivalent to the Riemann hypothesis.}\\


Here \cite{RobJMPA} is Robin's paper, which he published three years \emph{before} learning of  Ramanujan's Theorem.
However, a reader of \cite{RAM97} who neglects to look up \cite{RobJMPA} in the References is left with the misimpression ``that the above relation with $n_0=5040$ is equivalent to the Riemann hypothesis'' was proven \emph{after} whoever proved it learned of ``the above relation''!


%

\bigskip
 
In 1993, HCN Part~2 was submitted to Proceedings of the London Mathematical Society, which had published Part~1 in 1915. The paper was accepted, but could not be published, because Trinity College did not own the rights to Ramanujan's papers and was not able to obtain permission from his widow, Janaki.

\centerline{\includegraphics[width=1.5in]{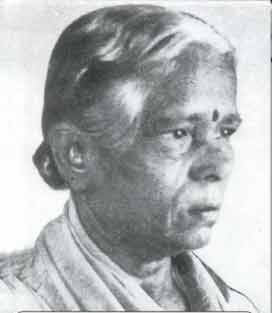}}
\begin{center}{ S. Janaki Ammal (1899--1994), Mrs. Ramanujan}\end{center}

\bigskip

Janaki passed away in 1994, and the paper was eventually published by Alladi in the first volume of his newly-founded Ramanujan Journal.\\

\centerline{\includegraphics[width=1.5in]{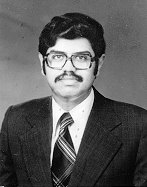}}\begin{center}{ Krishnaswami ALLADI, founder of The Ramanujan Journal}\end{center}

\bigskip

\centerline{\includegraphics[width=7.8in]{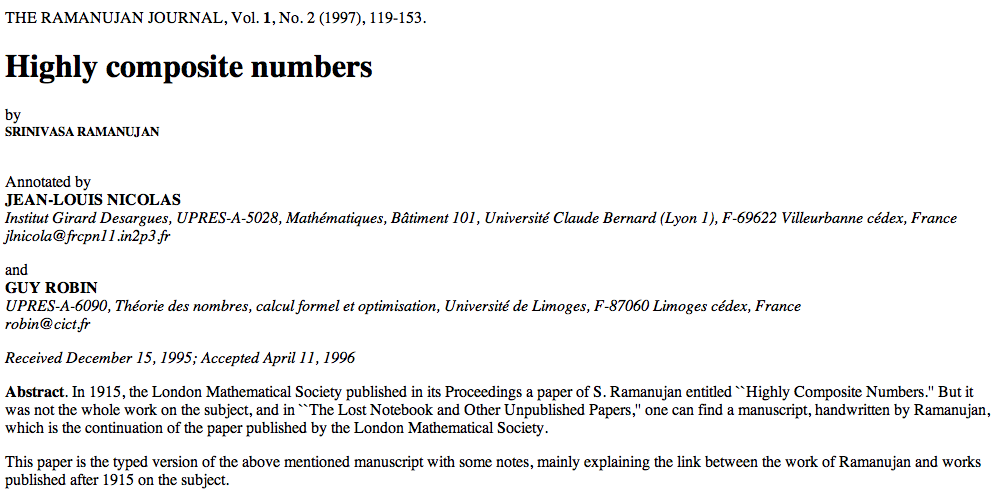}}

\begin{center}
{Ramanujan's HCN Part 2, annotated by Nicolas and Robin, The Ramanujan Journal, 1997}
\end{center}\


\smallskip
Here our story ends. If it has offended anyone, we apologize.


%
\

\vspace{2mm}
\vspace{2mm}


\centerline{\includegraphics[width=3in]{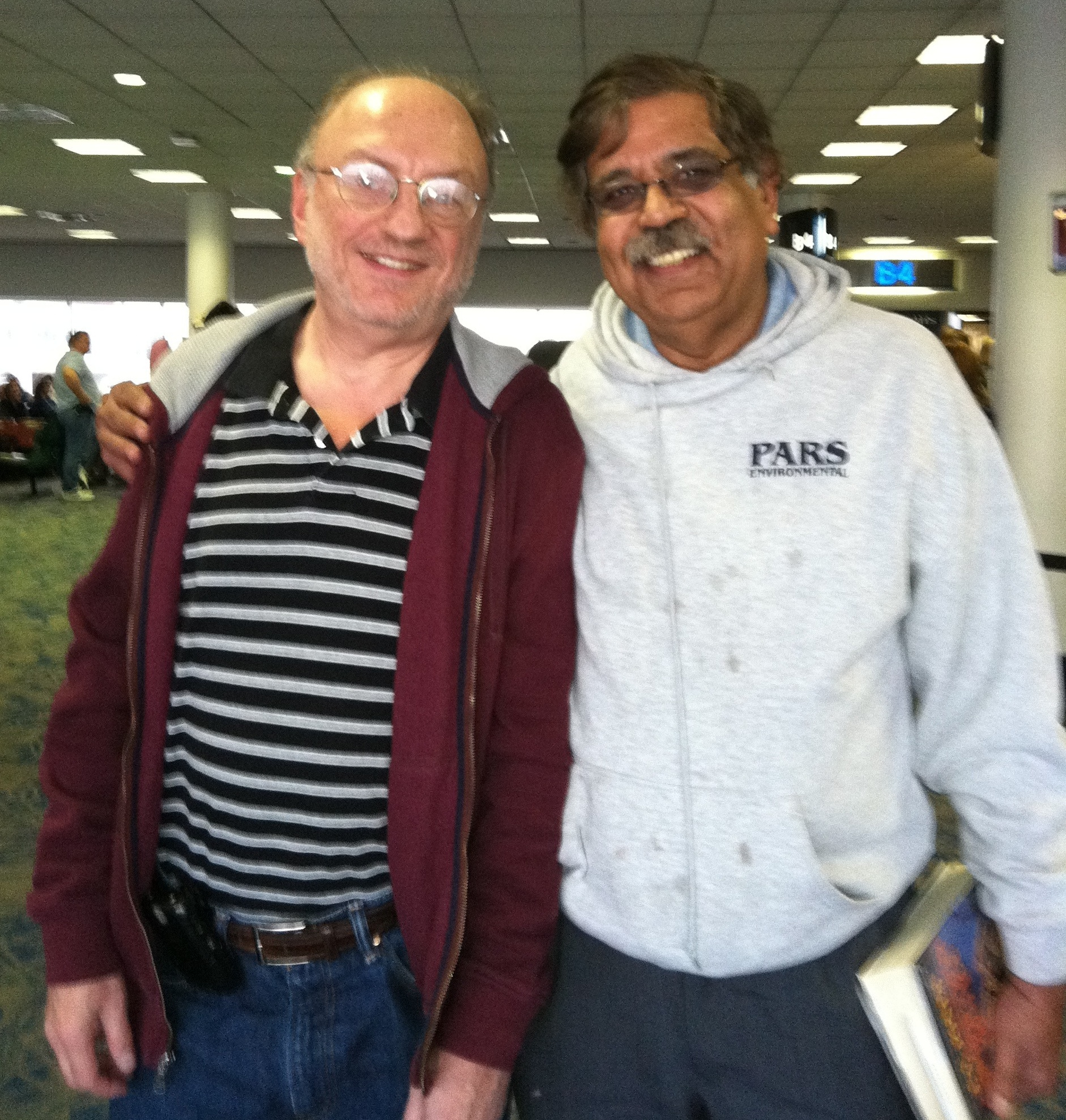}}
\begin{center}{Jonathan SONDOW (left) and Ramjee RAGHAVAN, Ramanujan's grandnephew, by chance (!) seatmates on a flight to Charlotte. (S.~flew on to RAMA125 in Gainesville, and R.~to Chicago.)}\end{center}

\bigskip

\def\refname{References}

\end{document}